\documentclass[11pt]{article}
\usepackage{epic,latexsym,amssymb}
\usepackage{amsfonts}
\usepackage{amscd}
\usepackage{amsmath}
\usepackage{graphicx}
\usepackage{color}
\usepackage{caption,
caption}
\usepackage{tikz}
\usetikzlibrary{arrows.meta}

\usepackage{float}

\usepackage{lineno}

\newtheorem{thm}{Theorem}

\newtheorem{lem}[thm]{Lemma}

\newtheorem{cor}[thm]{Corollary}

\newtheorem{claim}{Claim}

\newcommand{\cp}{\, \Box \,}

\newcommand{\proof}{\noindent\textbf{Proof. }}

\newcommand{\smallqed}{{\tiny ($\Box$)}}
\newcommand{\qed}{$\Box$}

\newcommand{\QEDmark}{\mbox{\textsc{qed}}}
\newcommand{\proofStarter}[1]{\textsc{#1} }

\def\vertex(#1){\put(#1){\circle*{2}}}
\def\vertexo(#1){\put(#1){\circle{2}}}
\def\vert(#1){\put(#1){\circle*{1.5}}}
\def\verto(#1){\put(#1){\circle{1.5}}}
\def\lab(#1)#2{\put(#1){\makebox(0,0)[c]{#2}}}
\definecolor{DarkGreen}{rgb}{0.2, 0.6, 0.3}

\definecolor{electricindigo}{rgb}{0.44, 0.0, 1.0}

\let\oldenumerate\enumerate
\renewcommand{\enumerate}{
  \oldenumerate
  \setlength{\itemsep}{0.5pt}
  \setlength{\parskip}{0pt}
  \setlength{\parsep}{0pt}
}

\makeatletter
\newcommand{\Rmnum}[1]{\expandafter\@slowromancap\romannumeral #1@}
\makeatother

\begin{document}

\title{Zero Forcing in Claw-Free Cubic Graphs}
\author{$^{1,2}$Randy Davila and $^{1}$Michael A. Henning\thanks{Research supported in part by the University of Johannesburg.} \\
\\
$^1$Department of Pure and Applied Mathematics\\
University of Johannesburg \\
Auckland Park 2006, South Africa \\
\small {\tt Email: mahenning@uj.ac.za} \\
\\
$^2$Department of Mathematics and Statistics\\
University of Houston--Downtown \\
Houston, TX 77002, USA \\
\small {\tt Email: davilar@uhd.edu}
}

\date{}
\maketitle

\begin{abstract}
The \emph{zero forcing number} of a simple graph, written $Z(G)$, is a NP-hard graph invariant which is the result of the \emph{zero forcing color change rule}. This graph invariant has been heavily studied by linear algebraists, physicists, and graph theorist. It's broad applicability and interesting combinatorial properties have attracted the attention of many researchers. Of particular interest, is that of bounding the zero forcing number from above. In this paper we show a surprising relation between the zero forcing number of a graph and the independence number of a graph, denoted $\alpha(G)$. Our main theorem states that if $G \ne K_4$ is a connected, cubic, claw-free graph, then $Z(G) \le \alpha(G) + 1$. This improves on best known upper bounds for $Z(G)$, as well as known lower bounds on $\alpha(G)$. As a consequence of this result, if $G \ne K_4$ is a connected, cubic, claw-free graph with order $n$, then $Z(G) \le \frac{2}{5}n + 1$. Additionally, under the hypothesis of our main theorem, we further show $Z(G) \le \alpha'(G)$, where $\alpha'(G)$ denotes the matching number of $G$.
\end{abstract}


\newpage
\section{Introduction}

Dynamic colorings in graphs, i.e., vertex (or edge) colorings that spread during discrete time intervals, have shown increasing relevance and applicability in the study of graphs. One of the most heavily studied dynamic colorings is that due to the \emph{zero forcing process}, and its associated graph invariant, the \emph{zero forcing number}. These notions were originally introduced in~\cite{AIM-Workshop} during a workshop on linear algebra, and have since, found relationships to well studied graph parameters such as the \emph{chromatic number},  the \emph{connected domination number}, the \emph{diameter}, and the \emph{independence number}, see for example~\cite{k-Forcing, DaHeTF1, Davila_Kenter, Chromatic}.

We next recall the zero forcing process as defined in~\cite{DaHeTF1}: Let $G$ be a finite and simple graph with vertex set $V(G)$, and let $S\subseteq V(G)$ be a set of initially ``colored'' vertices, while all other vertices being ``uncolored''. All vertices contained in $S$ are said to be $S$-colored while all vertices not in $S$ are $S$-uncolored. At each discrete time step, if a colored vertex has exactly one uncolored neighbor, then this colored vertex \emph{forces} its uncolored neighbor to become colored. If $v$ is such a colored vertex, then we call $v$ a forcing vertex, and say that $v$ has been \emph{played}. The initial set of colored vertices $S$ is a \emph{zero forcing set}, if by iteratively applying the above forcing rule, all of $V(G)$ becomes colored. We call such a set an $S$-\emph{forcing set}. If $S$ is a $S$-forcing set of $G$, and $v$ is a $S$-colored vertex which is played in the forcing process, then $v$ is a \emph{$S$-forcing vertex}. The \emph{zero forcing number}, written $Z(G)$, is the cardinality of a minimum zero forcing set in~$G$.

If $S$ is a forcing set in a graph $G$ and the subgraph $G[S]$ induced by $S$ contains no isolated vertex, then $S$ is a \emph{total forcing set}, abbreviated as a \emph{TF}-\emph{set} of $G$. The \emph{total forcing number} of $G$, written $F_t(G)$, is the cardinality of a minimum TF-set in $G$. The concept of a total forcing set was first introduced by Davila in~\cite{Davila}, and studied further, for example, in~\cite{DaHeTF1,DaHeTF2,DaHeTF3}.

In this paper, we study zero forcing sets in cubic, claw-free graphs. We proceed as follows. In Section~\ref{S:notation}, we give the necessary graph theory notation and terminology. Thereafter, we present our main results in Section~\ref{S:main}. In Section~\ref{S:known}, we present some known results. A proof of our main result is given in Section~\ref{S:proof1}.


\subsection{Notation and Terminology.} 
\label{S:notation}

For notation and terminology, we will typically follow~\cite{Henning}. Specifically, let $G$ be a graph with vertex set $V(G)$ and edge set $E(G)$. The order and size of $G$ will be denoted $n = |V(G)|$ and $m=|E(G)|$, respectively. A \emph{neighbor} of a vertex $v$ in $G$ is a vertex $u$ that is adjacent to $v$, that is, $uv \in E(G)$. The \emph{open neighborhood} of a vertex $v$ in $G$ is the set of neighbors of $v$, denoted $N_G(v)$. We denote the \emph{degree} of $v$ in $G$ by $d_G(v) = |N_G(v)|$. The minimum and maximum vertex degrees of $G$ will be denoted by $\delta(G)$ and $\Delta(G)$, respectively. A \emph{cubic graph} (also called a $3$-\emph{regular graph}) is a graph in which every vertex has degree~$3$.

Two edges in a graph $G$ are \emph{independent} if they are not adjacent in $G$. A set of pairwise independent edges of $G$ is called a \emph{matching} in $G$, while a matching of maximum cardinality is a \emph{maximum matching}. The number of edges in a maximum matching of $G$ is the \emph{matching number} of $G$, denoted $\alpha'(G)$. Matchings in graphs are extensively studied in the literature (see, for example, the classical book on matchings by Lov\'{a}sz and Plummer~\cite{LoPl86}, and the excellent survey articles by Plummer~\cite{Pl03} and Pulleyblank~\cite{Pu95}).

Two vertices in a graph $G$ are \emph{independent} if they are not neighbors. A set of pairwise independent vertices in $G$ is an \emph{independent set} of $G$. The number of vertices in a maximum independent set in $G$ is the \emph{independence number} of $G$, denoted $\alpha(G)$. We remark that the independence number is one of the most extensively studied graph invariants, see, for example \cite{Caro_one, Caro_two, Caro_Hansberg, Favaron}.

For a set of vertices $S\subseteq V(G)$, the subgraph induced by $S$ is denoted by $G[S]$. If $v\in V(G)$, we denote the graph obtained by deleting $v$ in $G$ by $G-v$. We denote the path, cycle, and complete graph on $n$ vertices by $P_n$, $C_n$, and $K_4$, respectively. A \emph{triangle} in $G$ is an induced subgraph of $G$ isomorphic to $K_3$, whereas a \emph{diamond} in $G$ is a subgraph of $G$ isomorphic to $K_4$ with one edge missing. A graph $G$ is \emph{$F$-free} if $G$ does not contain $F$ as an induced subgraph. In particular, if $G$ is $F$-free, where $F  = K_{1,3}$, then $G$ is \emph{claw-free}. Claw-free graphs are heavily studied and an excellent survey of claw-free graphs has been written by Flandrin, Faudree, and Ryjacek~\cite{claw_free_survey}. More recently, Chudnovsky and Seymour published a series of excellent papers in \emph{Journal of Combinatorial Theory Series B} on this topic~\cite{claw_free_JCTB}. 
We use the standard notation $[k] = \{1,2,\ldots,k\}$.

In this paper, we study zero forcing in connected, cubic, claw-free graphs. We proceed as follows. In Section~\ref{S:known}, we present some known results and a preliminary lemma. In Section~\ref{S:main}, we give our main result, namely Theorem~\ref{main_thm}. A proof of Theorem~\ref{main_thm} is given in Section~\ref{S:proof1}.

\section{Known Results and Preliminary Lemma} 
\label{S:known}

In this section, we present some known results and a preliminary lemma that will prove useful in proving our main result. Faudree et al.~\cite{claw_free_independence}, established the following upper bound on the independence number of a claw-free, cubic graph.

\begin{thm}{\rm (\cite{claw_free_independence})}
\label{alpha}
If $G$ is a claw-free, cubic graph of order $n$, then $\alpha(G) \le \frac{2}{5}n$.
\end{thm} 

Computation of the zero forcing number is known to be NP-hard~\cite{zf_np}, and as such, determining sharp upper and lower bounds on $Z(G)$ has attracted a considerable amount of interest. For example, Amos et al.~\cite{k-Forcing} showed that if $G$ is an isolate-free graph of order $n \ge 2$ and maximum degree $\Delta$, then $Z(G) \le ( \frac{\Delta}{\Delta+1} ) n$. Imposing the added restrictions that $G$ is connected and $\Delta \ge 2$, this bound is improved in~\cite{k-Forcing} to $Z(G) \le \frac{(\Delta -2)n+2}{\Delta -1}$. In the special case that $G$ is connected and cubic, this result simplifies to the following result.

\begin{thm}{\rm (\cite{k-Forcing})} 
\label{t:known1}
If $G$ is a connected, claw-free, cubic graph of order~$n$, then $Z(G) \le \frac{1}{2}n + 1$.
\end{thm}

It was shown in~\cite{DaHeTF3} that if $G \ne K_4$ is a connected, claw-free, cubic graph of order~$n$, then $F_t(G) \le \frac{1}{2}n$ and this bound is tight. Further, the (infinite family of) extremal graphs achieving equality in this bound are also characterized in~\cite{DaHeTF3}. As a consequence of this result, we have the following upper bound on the zero forcing number of a connected, claw-free, cubic graph. 

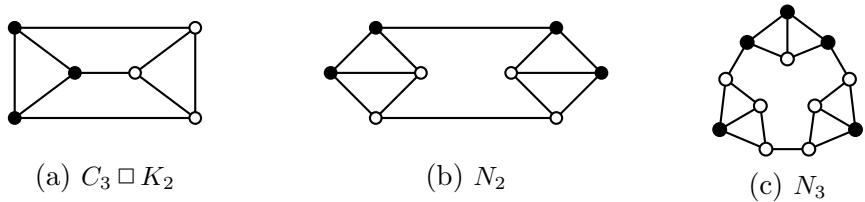
\begin{figure}[htb]
\begin{center}
\begin{tikzpicture}[scale=.8,style=thick,x=1cm,y=1cm]
\def\vr{2.75pt} 
\path (0,0) coordinate (b);
\path (0,1.5) coordinate (a);
\path (1,0.75) coordinate (c);
\path (2,0.75) coordinate (d);
\path (3,0) coordinate (f);
\path (3,1.5) coordinate (e);
\path (5.25,0.75) coordinate (g);
\path (6.75,0.75) coordinate (j);
\path (8.25,0.75) coordinate (k);
\path (9.75,0.75) coordinate (n);
\path (6,0) coordinate (i);
\path (9,0) coordinate (m);
\path (6,1.5) coordinate (h);
\path (9,1.5) coordinate (l);
\draw (a) -- (b);
\draw (a) -- (c);
\draw (a) -- (e);
\draw (b) -- (c);
\draw (b) -- (f);
\draw (c) -- (d);
\draw (d) -- (e);
\draw (d) -- (f);
\draw (e) -- (f);
\draw (g) -- (h);
\draw (g) -- (i);
\draw (g) -- (j);
\draw (h) -- (j);
\draw (h) -- (l);
\draw (i) -- (j);
\draw (i) -- (m);
\draw (k) -- (l);
\draw (k) -- (m);
\draw (k) -- (n);
\draw (l) -- (n);
\draw (m) -- (n);
\draw (a) [fill=black] circle (\vr);
\draw (b) [fill=black] circle (\vr);
\draw (c) [fill=black] circle (\vr);
\draw (d) [fill=white] circle (\vr);
\draw (e) [fill=white] circle (\vr);
\draw (f) [fill=white] circle (\vr);
\draw (g) [fill=black] circle (\vr);
\draw (h) [fill=black] circle (\vr);
\draw (i) [fill=white] circle (\vr);
\draw (j) [fill=white] circle (\vr);
\draw (k) [fill=white] circle (\vr);
\draw (l) [fill=black] circle (\vr);
\draw (m) [fill=white] circle (\vr);
\draw (n) [fill=black] circle (\vr);
\draw (1.5,-1) node {(a) {\small $C_3 \cp K_2$}};
\draw (7.5,-1) node {(b) {\small $N_2$}};
\end{tikzpicture}
\hspace*{1cm}
\tikzstyle{every node}=[circle, draw, fill=black!0, inner sep=0pt,minimum width=.17cm]
\begin{tikzpicture}[thick,scale=.5]
  \draw(0,0) { 
    +(1.80,4.65) -- +(1.80,3.40)
    +(1.80,3.40) -- +(2.88,3.84)
    +(2.88,3.84) -- +(1.80,4.65)
    +(1.80,4.65) -- +(0.73,3.84)
    +(0.73,3.84) -- +(1.80,3.40)
    +(0.73,3.84) -- +(0.16,2.86)
    +(0.16,2.86) -- +(0.00,1.52)
    +(0.00,1.52) -- +(1.08,2.15)
    +(1.08,2.15) -- +(0.16,2.86)
    +(1.08,2.15) -- +(1.23,1.00)
    +(1.23,1.00) -- +(0.00,1.52)
    +(1.23,1.00) -- +(2.37,1.00)
    +(2.37,1.00) -- +(2.53,2.15)
    +(2.53,2.15) -- +(3.61,1.52)
    +(3.61,1.52) -- +(2.37,1.00)
    +(2.53,2.15) -- +(3.45,2.86)
    +(3.45,2.86) -- +(3.61,1.52)
    +(3.45,2.86) -- +(2.88,3.84)
    +(1.80,4.65) node[fill=black]{}
    +(0.73,3.84) node[fill=black]{}
    +(2.88,3.84) node[fill=black]{}
    +(1.80,3.40) node{}
    +(0.16,2.86) node{}
    +(1.08,2.15) node{}
    +(1.23,1.00) node{}
    +(2.37,1.00) node{}
    +(2.53,2.15) node{}
    +(3.45,2.86) node{}
    +(0.00,1.52) node[fill=black]{}
    +(3.61,1.52) node[fill=black]{}
    +(1.80,0) node[rectangle, draw=white!0, fill=white!100]{{(c) {\small $N_3$}}}  
  };
\end{tikzpicture}
\end{center}
\vskip -0.65cm
\caption{The prism $C_3 \cp K_2$ and the diamond-necklaces $N_2$ and $N_3$.} \label{f:prism}
\end{figure}

\begin{thm}{\rm (\cite{DaHeTF3})} 
\label{t:known2}
If $G \ne K_4$ is a connected, claw-free, cubic graph of order $n$, then $Z(G) \le \frac{1}{2}n$ with equality if and only if $G$ is the prism $C_3 \cp K_2$ (shown in Figure~\ref{f:prism}(a)) or $G$ is the diamond-necklace $N_2$ (shown in Figure~\ref{f:prism}(b)). 
\end{thm}

We note that $Z(C_3 \cp K_2) = 3$, $Z(N_2) = 4$ and $Z(N_3) = 5$. Moreover, the darkened vertices shown in Figure~\ref{f:prism}(a),~\ref{f:prism}(b) and~\ref{f:prism}(c) form a minimum zero forcing set in the associated graph.

The following property of connected, claw-free, cubic graphs is established in~\cite{diamond_lem}.

\begin{lem}{\rm (\cite{diamond_lem})}
 \label{l:diamond_lem}
If $G \ne K_4$ is a connected, claw-free, cubic graph of order $n$, then the vertex set $V(G)$ can be uniquely partitioned into sets each of which induces a triangle or a diamond in $G$.
\end{lem}

By Lemma~\ref{l:diamond_lem}, the vertex set $V(G)$ of connected, claw-free, cubic graph $G \ne K_4$  can be uniquely partitioned into sets each of which induce a triangle or a diamond in $G$. Following the notation introduced in~\cite{diamond_lem}, we refer to such a partition as a \emph{triangle}-\emph{diamond partition} of $G$, abbreviated $\Delta$-D-partition. We call every triangle and diamond induced by a set in our $\Delta$-D-partition a \emph{unit} of the partition. A unit that is a triangle is called a \emph{triangle-unit} and a unit that is a diamond is called a \emph{diamond-unit}. (We note that a triangle-unit is a triangle that does not belong to a diamond.) We say that two units in the $\Delta$-D-partition are \emph{adjacent} if there is an edge joining a vertex in one unit to a vertex in the other unit.

\section{Main Result}
\label{S:main}

We have two immediate aims. First to establish a relationship between the zero forcing number of a cubic, claw-free graph and its independence and matching numbers. Secondly, to obtain a tight upper bound on the zero forcing number of a cubic, claw-free graph in terms of its order. More precisely, we shall prove the following results. A proof of Theorem~\ref{main_thm} is given in Section~\ref{S:proof1}.

\begin{thm}
\label{main_thm}
If $G \ne K_4$ is a connected, claw-free, cubic graph, then the following holds. \\
\indent {\rm (a)} $Z(G) \le \alpha(G)+1$.
\\
\indent {\rm (b)} $Z(G) \le \alpha'(G)$. 
\end{thm} 

We note that if $G$ is the prism $C_3 \cp K_2$ (shown in Figure~\ref{f:prism}(a)) or the diamond-necklace $N_2$ (shown in Figure~\ref{f:prism}(b)), then $Z(G) = \alpha'(G) = \alpha(G)+1$. Thus, the bounds of Theorem~\ref{main_thm} are achievable. If $G$ is the diamond-necklace $N_3$ (shown in Figure~\ref{f:prism}(c)), then $Z(G) = 5 = \alpha(G)+1$. As an immediate consequence Theorem~\ref{alpha} and Theorem~\ref{main_thm}(a), we have the following upper bound on the zero forcing number of a claw-free, cubic graph in terms of its order.

\begin{cor}
\label{c:main}
If $G \ne K_4$ is a connected, claw-free, cubic graph of order~$n$, then 
\[
Z(G) \le \frac{2}{5}n + 1.
\] 
\end{cor}

\section{Proof of Theorem~\ref{main_thm}}
\label{S:proof1}

In this section, we prove Theorem~\ref{main_thm}. Recall its statement. 

\medskip
\noindent \textbf{Theorem~\ref{main_thm}}. \emph{If $G \ne K_4$ is a connected, claw-free, cubic graph, then the following holds. \\
\indent {\rm (a)} $Z(G) \le \alpha(G)+1$.
\\
\indent {\rm (b)} $Z(G) \le \alpha'(G)$.}

\proof Let $G$ be a connected, claw-free, cubic graph of order~$n \ge 6$. If $n = 6$, then $G$ is the prism $C_3 \cp K_2$ and $Z(G) = \alpha'(G) = 3$ and $\alpha(G) = 2$. If $n = 8$, then $G$ is the diamond-necklace $N_2$ and $Z(G) = \alpha'(G) = 4$ and $\alpha(G) = 3$. Hence, we may assume in what follows that $n \ge 10$. We now consider the (unique) $\Delta$-D-partition of $G$ given by Lemma~\ref{l:diamond_lem}. We will greedily construct a zero forcing set, and while doing so we also produce an independent set of vertices. We remark that our technique relies on greedily coloring vertices which are independent of all but at most one vertex which has been previously greedily colored. Moreover, we also ensure that each greedily colored vertex is played during the forcing process on~$G$. We start this process with the following initialization which gives a set of colored vertices from which we start our greedy coloring process.

\begin{figure}[htb]
\begin{center}
\begin{tikzpicture}[scale=.8,style=thick,x=1cm,y=1cm, =>stealth]
\def\vr{2.5pt} 
%

\path (-2,0) coordinate (x1);
\path (0,1) coordinate (x2);
\path (0,-1) coordinate (x3);
\path (2,0) coordinate (x4);

\path (3,0) coordinate (y1);
\path (-3,0) coordinate (w1);

%


\draw (x2) -- (x3);
\draw (w1) -- (x1);
\draw (x1) -- (x2);
\draw (x2) -- (x4);
\draw (x4) -- (y1);
\draw (x1) -- (x3);
\draw (x2) -- (x3);
\draw (x4) -- (x3);

\draw (x1) [fill=black]circle (\vr);
\draw (x2) [fill=black] circle (\vr);
\draw (x3) [fill=white] circle (\vr);
\draw (x4) [fill=black] circle (\vr);

\draw (y1) [fill=white] circle (\vr);
\draw (w1) [fill=white]circle (\vr);

\draw[anchor = south] (w1) node {{\small $w_1$}};
\draw[anchor = south] (x1) node {{\small $x_1$}};

\draw[anchor = south] (x2) node {{\small $x_2$}};
\draw[anchor = south] (x4) node {{\small $x_4$}};
\draw[anchor = south] (y1) node {{\small $y_1$}};
\draw[anchor = north] (x3) node {{\small $x_3$}};
\draw (0,-2) node {(a)};
%

\path (6,0) coordinate (x1);
\path (8,1) coordinate (x2);
\path (8,-1) coordinate (x3);
\path (10,0) coordinate (x4);

\path (11,0) coordinate (y1);
\path (5,0) coordinate (w1);

%

\draw[black, arrows={->[line width=1pt,black,length=3mm,width=3mm]}] (x2) -- (x3);
\draw[black, arrows={->[line width=1pt,black,length=3mm,width=3mm]}] (x1) -- (w1);
\draw[black, arrows={->[line width=1pt,black,length=3mm,width=3mm]}] (x4) -- (y1);

\draw (w1) -- (x1);
\draw (x1) -- (x2);
\draw (x2) -- (x4);
\draw (x4) -- (y1);
\draw (x1) -- (x3);
\draw (x2) -- (x3);
\draw (x4) -- (x3);

\draw (x1) [fill=black]circle (\vr);
\draw (x2) [fill=black] circle (\vr);
\draw (x3) [fill=black] circle (\vr);
\draw (x4) [fill=black] circle (\vr);

\draw (y1) [fill=white] circle (\vr);
\draw (w1) [fill=white]circle (\vr);

\draw[anchor = south] (w1) node {{\small $w_1$}};
\draw[anchor = south] (x1) node {{\small $x_1$}};

\draw[anchor = south] (x2) node {{\small $x_2$}};
\draw[anchor = south] (x4) node {{\small $x_4$}};
\draw[anchor = south] (y1) node {{\small $y_1$}};
\draw[anchor = north] (x3) node {{\small $x_3$}};
\draw (8,-2) node {(b)};
\end{tikzpicture}
\end{center}
\vskip -0.5 cm
\caption{Initialization with the diamond-unit $D$.} \label{f:fig1}
\end{figure}
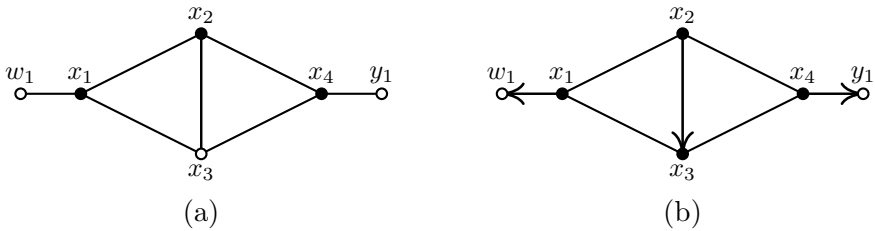

\noindent
\textbf{Initialize.} If $G$ contains a diamond-unit, we initialize as follows: Let $D$ be an arbitrary diamond-unit in $G$, where $V(D) = \{x_1, x_2, x_3, x_4\}$ and where $x_1x_4$ is the missing edge in $D$. Let $S = \{x_1, x_2, x_4\}$ be an initial set of colored vertices. We note that $x_1$ and $x_4$ have exactly one neighbor outside of $D$. Let $w_1$ and $y_1$ be the neighbors of $x_1$ and $x_4$, respectively, outside of $D$. By the claw-freeness of $G$, these neighbors are distinct. Under the coloring $S$, observe that $x_2$ may force $x_3$ to become colored. Allowing $x_2$ to force $x_3$ to become colored, we next observe that each of $x_1$ and $x_4$ has exactly one uncolored neighbor, namely, $w_1$ and $y_1$, respectively. Let $x_1$ and $x_4$ be played, and observe that all vertices in $D$ have become colored, along with one vertex from each unit adjacent to $D$. Moreover, $I = \{x_1, x_4\}\subseteq S$ forms an independent set, and each vertex from $I$ has been played. See Figure~\ref{f:fig1}(a) and~\ref{f:fig1}(b) for an illustration.

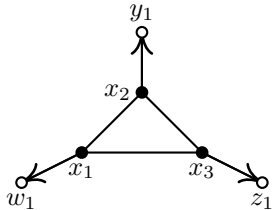
\begin{figure}[htb]
\begin{center}
\begin{tikzpicture}[scale=.8,style=thick,x=1cm,y=1cm, =>stealth]
\def\vr{2.5pt} 
%

\path (-1,0) coordinate (x1);
\path (0,1) coordinate (x2);
\path (1,0) coordinate (x3);

\path (2,-.5) coordinate (z1);
\path (0,2) coordinate (y1);
\path (-2,-.5) coordinate (w1);

%

\draw[black, arrows={->[line width=1pt,black,length=3mm,width=3mm]}] (x2) -- (y1);
\draw[black, arrows={->[line width=1pt,black,length=3mm,width=3mm]}] (x1) -- (w1);
\draw[black, arrows={->[line width=1pt,black,length=3mm,width=3mm]}] (x3) -- (z1);

\draw (w1) -- (x1);
\draw (x1) -- (x2);
\draw (x2) -- (y1);

\draw (x1) -- (x3);
\draw (x2) -- (x3);
\draw (z1) -- (x3);

\draw (x1) [fill=black]circle (\vr);
\draw (x2) [fill=black] circle (\vr);
\draw (x3) [fill=black] circle (\vr);

\draw (z1) [fill=white] circle (\vr);
\draw (y1) [fill=white] circle (\vr);
\draw (w1) [fill=white]circle (\vr);

\draw[anchor = north] (w1) node {{\small $w_1$}};
\draw[anchor = north] (x1) node {{\small $x_1$}};

\draw[anchor = east] (x2) node {{\small $x_2$}};
\draw[anchor = north] (z1) node {{\small $z_1$}};
\draw[anchor = south] (y1) node {{\small $y_1$}};
\draw[anchor = north] (x3) node {{\small $x_3$}};

%



%

\end{tikzpicture}
\end{center}
\vskip -0.5 cm
\caption{Initialization with the triangle-unit $T$.} \label{f:fig2}
\end{figure}

If $G$ does not contain a diamond-unit, we initialize as follows: Let $T$ be an arbitrary triangle-unit in $G$, where $V(T) = \{x_1, x_2, x_3\}$. Since $G$ does not contain a diamond-unit, we note that every vertex of $G$ is contained in a triangle-unit. In particular, this implies that no two vertices in $T$ have a common neighbor outside of $T$. Let $w_1$, $y_1$, and $z_1$ be the neighbors of $x_1$, $x_2$, and $x_3$, respectively, outside of $T$. Let $S = \{x_1, x_2, x_3\}$ be a set of initially colored vertices. Under the coloring $S$, observe that each $S$-colored vertex has exactly one $S$-uncolored neighbor. Let each vertex in $S$ force their respective $S$-uncolored neighbor, i.e., allow $w_1$, $y_1$, and $z_1$ to become colored. If $w_1$, $y_1$, and $z_1$ all belong to the same triangle-unit, then $G$ is the prism $C_3 \cp K_2$, contradicting our assumption that $n \ge 10$. Renaming vertices if necessary, we may therefore assume without loss of generality that $w_1$ and $y_1$ lie in distinct units. In this case, we note that $I = \{x_1, y_1\}$ forms an independent set, where $x_1$ is a played vertex during the forcing process on $G$. See Figure~\ref{f:fig2} for an illustration.

\newpage
In each case of the initialization step, we start with a set of colored vertices $S$, where all vertices are contained in a common unit, and each vertex forces a vertex in a neighboring unit. We call this unit containing the vertices of $S$ an \emph{intially-forcing unit}. Note that if the initially-forcing unit is a diamond-unit, then all vertices of $I$ (currently) are played during the forcing process on $G$, and if the initially-forcing unit is a triangle-unit, then all except possibly one vertex in $I$, namely the vertex $y_1$, is played during the forcing process on $G$. Our next step is to greedily add vertices to both $S$ and $I$. Moreover, along the way we will specify exactly how our greedily colored vertices will be allowed to force during the forcing process on $G$.

\noindent
\textbf{Greedy Coloring.} Let $S$ be defined as in the initialization step, and let $U_1$ denote the initially-forcing unit. We adopt our earlier notation as defined in the initialization section. Our key notion is that we add vertices to $S$ and $I$, so that we may ensure each vertex contained in $I$ is also contained in $S$, and further, that each vertex in $I$ may be played during the forcing process on $G$.

\begin{claim}
\label{one}
If $U_1$ is a diamond-unit, then $Z(G) \le \alpha(G) + 1$ and $Z(G) \le \alpha'(G)$. 
\end{claim}
\proof Suppose that $U_1$ is a diamond-unit. Adopting our earlier notation, recall that $I = \{x_1, x_4\}$ and $S = \{x_1, x_2, x_4\}$. Starting from the set $S \subseteq V(U_1)$, let the forcing process propagate throughout $V(G)$. If all of $V(G)$ becomes colored, then we are done since in this case $Z(G) \le |S| = 3 = |I| + 1 \le \alpha(G)$. Hence we may assume that starting with the set $S$, the forcing process halts before all of $V(G)$ becomes colored, for otherwise the claim is satisfied. This implies that at some point of the forcing process, no further forcing steps will occur. Thus, there must be a colored vertex, say $v$, with exactly two uncolored neighbors. Note that so far, each vertex in $I = \{x_1, x_4\}$ has been played, and so neither uncolored neighbor of $v$ belongs to $I$ or is adjacent to a vertex of $I$. 
Moreover, we will assume that at this initial stage of the forcing process, we have not greedily colored any vertices. We now apply the following rules where the only vertices colored by our process (we exclude vertices colored by the forcing process) are the vertices contained in $U_1$.

\begin{figure}[htb]
\begin{center}
\begin{tikzpicture}[scale=.8,style=thick,x=1cm,y=1cm, =>stealth]
\def\vr{2.5pt} 
%

\path (-2,0) coordinate (v);
\path (-1,1) coordinate (w);
\path (0,0) coordinate (y);

\path (1,-.5) coordinate (z1);
\path (-1,2) coordinate (y1);
\path (-3,-.5) coordinate (w1);

%

\draw[black, arrows={->[line width=1pt,black,length=3mm,width=3mm]}] (w1) -- (v);

\draw (w1) -- (v);
\draw (v) -- (w);
\draw (w) -- (y1);

\draw (v) -- (y);
\draw (w) -- (y);
\draw (z1) -- (y);

\draw (v) [fill=black]circle (\vr);
\draw (w) [fill=white] circle (\vr);
\draw (y) [fill=white] circle (\vr);

\draw (z1) [fill=white] circle (\vr);
\draw (y1) [fill=white] circle (\vr);
\draw (w1) [fill=black]circle (\vr);

\draw (-1,0.4) node {{\small $T_v$}};
\draw[anchor = north] (v) node {{\small $v$}};
\draw[anchor = east] (w) node {{\small $w$}};
\draw[anchor = west] (y1) node {{\small $w'$}};
\draw[anchor = north] (y) node {{\small $y$}};

%



%


\path (3,0) coordinate (v);
\path (4,1) coordinate (w);
\path (5,0) coordinate (y);

\path (6,-.5) coordinate (z1);
\path (4,2) coordinate (y1);
\path (2,-.5) coordinate (w1);

%

\draw[black, arrows={->[line width=1pt,black,length=3mm,width=3mm]}] (w1) -- (v);
\draw[black, arrows={->[line width=1pt,black,length=3mm,width=3mm]}] (w) -- (y);

\draw (w1) -- (v);
\draw (v) -- (w);
\draw (w) -- (y1);

\draw (v) -- (y);
\draw (w) -- (y);
\draw (z1) -- (y);

\draw (v) [fill=black]circle (\vr);
\draw (w) [fill=red] circle (\vr);
\draw (y) [fill=white] circle (\vr);

\draw (z1) [fill=white] circle (\vr);
\draw (y1) [fill=black] circle (\vr);
\draw (w1) [fill=black]circle (\vr);

\draw (-1,-1.25) node {(a)};
\draw (4,-1.25) node {(b)};
\draw[anchor = north] (v) node {{\small $v$}};
\draw[anchor = west] (y1) node {{\small $w'$}};
\draw[anchor = east] (w) node {{\small $w$}};
\draw[anchor = north] (y) node {{\small $y$}};

%



%


\path (8,0) coordinate (v);
\path (9,1) coordinate (w);
\path (10,0) coordinate (y);

\path (11,-.5) coordinate (z1);
\path (9,2) coordinate (y1);
\path (7,-.5) coordinate (w1);

%

\draw[black, arrows={->[line width=1pt,black,length=3mm,width=3mm]}] (w) -- (y1);
\draw[black, arrows={->[line width=1pt,black,length=3mm,width=3mm]}] (w1) -- (v);
\draw[black, arrows={->[line width=1pt,black,length=3mm,width=3mm]}] (v) -- (y);

\draw (w1) -- (v);
\draw (v) -- (w);
\draw (w) -- (y1);

\draw (v) -- (y);
\draw (w) -- (y);
\draw (z1) -- (y);

\draw (v) [fill=black]circle (\vr);
\draw (w) [fill=red] circle (\vr);
\draw (y) [fill=black] circle (\vr);

\draw (z1) [fill=white] circle (\vr);
\draw (y1) [fill=white] circle (\vr);
\draw (w1) [fill=black]circle (\vr);

\draw (9,-1.25) node {(c)};
\draw[anchor = north] (v) node {{\small $v$}};
\draw[anchor = west] (y1) node {{\small $w'$}};
\draw[anchor = east] (w) node {{\small $w$}};
\draw[anchor = north] (y) node {{\small $y$}};

\end{tikzpicture}
\end{center}
\vskip -0.5 cm
\caption{Illustration of the Triangle-Rule applied to $T_v$.} \label{Triangle_Rule}
\end{figure}
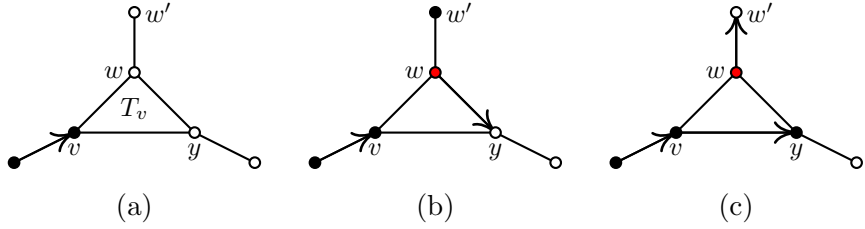

\noindent
\textbf{Triangle-Rule.} Suppose that $v$ is contained in a triangle-unit, say $T_v$ where $V(T_v) = \{v,w,y\}$. By our earlier assumptions, both $w$ and $y$ are currently uncolored. Moreover, neither $w$ nor $y$ are adjacent to any vertex in $I$, since $U_1$ being a diamond-unit implies that currently all vertices of $I$ are colored and have been played. We now greedily color the vertex $w$, and update $S$ by adding to it the vertex~$w$; that is, $S:=S\cup\{w\}$. Further, we also update $I$ by adding to it the vertex $w$; that is, $I := I\cup\{w\}$. Let $w'$ be the neighbor of $w$ not in $T_v$. If $w'$ is a colored vertex, then the vertex $w$ may force $y$ to become colored, as illustrated in Figure~\ref{Triangle_Rule}(b). If $w'$ is not a colored vertex, then the vertex $v$ may force $y$ to become colored, and thereafter the vertex $w$ may force the vertex $w'$ to become colored, as illustrated in Figure~\ref{Triangle_Rule}(c). In both cases, $w$ is a greedily colored vertex (which is colored red in Figure~\ref{Triangle_Rule}) that is played during the forcing process on $G$, and the updated set $I$ remains an independent set.

\begin{figure}[htb]
\begin{center}
\begin{tikzpicture}[scale=.8,style=thick,x=1cm,y=1cm, =>stealth]
\def\vr{2.5pt} 

\path (-2,1) coordinate (v);
\path (-1,2) coordinate (w);
\path (-1,0) coordinate (y);
\path (0,1) coordinate (z);

\path (-3, 1) coordinate (x1);
\path(1,1) coordinate (x2);

%

\draw[black, arrows={->[line width=1pt,black,length=3mm,width=3mm]}] (x1) -- (v);

\draw (v) -- (w);
\draw (w) -- (z);
\draw (v) -- (y);
\draw (w) -- (y);
\draw (z) -- (y);

\draw (x1)--(v);
\draw(x2)--(z);

\draw (v) [fill=black]circle (\vr);
\draw (w) [fill=white] circle (\vr);
\draw (y) [fill=white] circle (\vr);
\draw (z) [fill=white] circle (\vr);

\draw(x1) [fill=black]circle (\vr);
\draw(x2) [fill=white]circle (\vr);

\draw (-1,-1) node {(a)};

\draw[anchor = north] (v) node {{\small $v$}};

\draw[anchor = east] (w) node {{\small $w$}};
\draw[anchor = north] (z) node {{\small $z$}};
\draw[anchor = north] (y) node {{\small $y$}};


\path (3,1) coordinate (v);
\path (4,2) coordinate (w);
\path (4,0) coordinate (y);
\path (5,1) coordinate (z);

\path (2, 1) coordinate (x1);
\path(6,1) coordinate (x2);

%

\draw[black, arrows={->[line width=1pt,black,length=3mm,width=3mm]}] (x1) -- (v);
\draw[black, arrows={->[line width=1pt,black,length=3mm,width=3mm]}] (w) -- (y);

\draw (v) -- (w);
\draw (w) -- (z);
\draw (v) -- (y);
\draw (w) -- (y);
\draw (z) -- (y);

\draw (x1)--(v);
\draw(x2)--(z);

\draw (v) [fill=black]circle (\vr);
\draw (w) [fill=red] circle (\vr);
\draw (y) [fill=white] circle (\vr);
\draw (z) [fill=black] circle (\vr);

\draw(x1) [fill=black]circle (\vr);
\draw(x2) [fill=white]circle (\vr);

\draw (4,-1) node {(b)};
\draw[anchor = north] (v) node {{\small $v$}};

\draw[anchor = east] (w) node {{\small $w$}};
\draw[anchor = north] (z) node {{\small $z$}};
\draw[anchor = north] (y) node {{\small $y$}};

\path (8,1) coordinate (v);
\path (9,2) coordinate (w);
\path (9,0) coordinate (y);
\path (10,1) coordinate (z);

\path (7, 1) coordinate (x1);
\path(11,1) coordinate (x2);

%

\draw[black, arrows={->[line width=1pt,black,length=3mm,width=3mm]}] (v) -- (y);
\draw[black, arrows={->[line width=1pt,black,length=3mm,width=3mm]}] (w) -- (z);
\draw[black, arrows={->[line width=1pt,black,length=3mm,width=3mm]}] (x1) -- (v);

\draw (v) -- (w);
\draw (w) -- (z);
\draw (v) -- (y);
\draw (w) -- (y);
\draw (z) -- (y);

\draw (x1)--(v);
\draw(x2)--(z);

\draw (v) [fill=black]circle (\vr);
\draw (w) [fill=red] circle (\vr);
\draw (y) [fill=black] circle (\vr);
\draw (z) [fill=white] circle (\vr);

\draw(x1) [fill=black]circle (\vr);
\draw(x2) [fill=white]circle (\vr);

\draw (9,-1) node {(c)};
\draw[anchor = north] (v) node {{\small $v$}};

\draw[anchor = east] (w) node {{\small $w$}};
\draw[anchor = north] (z) node {{\small $z$}};
\draw[anchor = north] (y) node {{\small $y$}};

\end{tikzpicture}
\end{center}
\vskip -0.5 cm
\caption{Illustration of the Diamond-Rule applied to $D_v$.} \label{Diamond_Rule}
\end{figure}
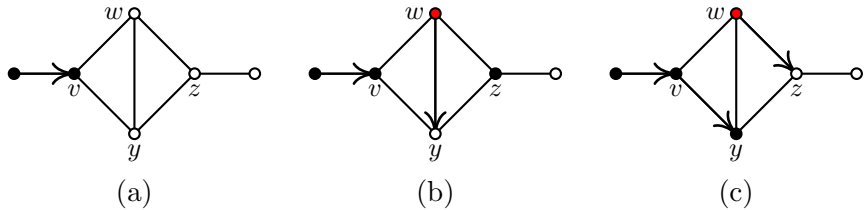

\noindent
\textbf{Diamond-Rule.} Next suppose $v$ is contained in a diamond-unit, say $D_v$ where $V(D_v) = \{v,w,y, z\}$ and $vz$ is the missing edge in $D_v$ as illustrated in Figure~\ref{Diamond_Rule}(a). Since $v$ is colored, and the forcing process halted at $v$, we note that both $w$ and $y$ are currently uncolored. Thus far, each of our initially colored vertices in $S$ (and $I$) have been played, and so, we are assured that $D_v$ contains no vertices from $S$. Since no vertices in $D_v$ are contained in $S$, we observe that both $w$ and $y$ are independent from vertices in $I$. We now greedily color the vertex $w$, and update $S$ by adding to it the vertex~$w$; that is, $S:=S\cup\{w\}$. Further, we also update $I$ by adding to it the vertex $w$; that is, $I := I\cup\{w\}$. If $z$ is colored, then $w$ may force $y$ to become colored, as illustrated in  Figure~\ref{Diamond_Rule}(b). Otherwise, if $z$ is not colored, then the vertex $v$ may first force $y$ to become colored, and thereafter the vertex $w$ may force $z$ to become colored, as illustrated in  Figure~\ref{Diamond_Rule}(c). In both cases, $w$ is a greedily colored vertex (which is colored red in Figure~\ref{Diamond_Rule}) that is played during the forcing process on $G$, and the updated set $I$ remains an independent set.

We now allow the zero forcing process to continue, and at each halting point (before all of $V(G)$ becomes colored) there must be a colored vertex, say $v$, with exactly two uncolored neighbors. We apply either the Triangle-Rule, or the Diamond-Rule, according to whether the vertex $v$ belongs to a triangle-unit or a diamond-unit, respectively. Since $G$ is connected, and since each rule allows the forcing process to continue, we are assured that our greedy coloring process will result in all of $V(G)$ becoming colored. Moreover, both the Triangle-Rule and the Diamond-Rule ensure that at each halting point of the forcing process, we greedily color a vertex which is independent from any previously colored vertices. Indeed, the only vertex which we have colored that is not independent from other colored vertices is the vertex $x_2$ from the initially-forcing unit $U_1$. It follows that $I$ is an independent set and $S$ is a zero forcing set with $|S| = |I|+1$, implying that $Z(T) \le |S| = |I| + 1 \le \alpha(G) + 1$. This completes the proof of part~(a).

To prove part~(b), we note that every vertex in the constructed set $S$ is an $S$-forcing vertex. Since a vertex in $S$ only forces one new vertex to be colored, and since no two distinct vertices in $S$ force the same vertex to be colored, the edges along which the vertices in $S$ force are independent of each other, implying that the graph $G$ contains a matching of size~$|S|$. Hence, $Z(T) \le |S| \le \alpha'(G)$. This completes the proof of Claim~\ref{one}.~\smallqed

\medskip
If the graph $G$ contains a diamond-unit, then our initialization process would have chosen a diamond-unit as the initially-forcing unit $U_1$, and the desired result would follow by Claim~\ref{one}. Hence, we may assume that $G$ contains no diamond-unit, for otherwise there is nothing left to prove. Thus, our initial-forcing unit $U_1$ is a triangle-unit. Adopting our earlier notation, recall that $V(U_1) = \{x_1, x_2, x_3\}$, where $w_1$, $y_1$, and $z_1$ are the neighbors of $x_1$, $x_2$, and $x_3$, respectively, outside of $U_1$. Let $U_2$,  $U_3$ and $U_4$ be the triangle-units containing $w_1$, $y_1$ and $z_1$, respectively. By our earlier assumptions, the units $U_2$ and $U_3$ are distinct units.

\begin{figure}[htb]
\begin{center}
\begin{tikzpicture}[scale=.8,style=thick,x=1cm,y=1cm, =>stealth]
\def\vr{2.5pt} 

\path (3,0) coordinate (x1);
\path (4,1) coordinate (x2);
\path (5,0) coordinate (x3);

\path (-2,0) coordinate (v1);
\path (-1,1) coordinate (v2);
\path (0,0) coordinate (v3);

\path (8,0) coordinate (y2);
\path (9,1) coordinate (y1);
\path (10,0) coordinate (y3);


\draw (x1) -- (x2);
\draw (x1) -- (x3);
\draw (x2) -- (x3);

\draw (v1) -- (v2);
\draw (v1) -- (v3);
\draw (v2) -- (v3);

\draw(y1)--(y2);
\draw(y3)--(y2);

\draw[black, arrows={->[line width=1pt,black,length=3mm,width=3mm]}] (x2) -- (y1);
\draw[black, arrows={->[line width=1pt,black,length=3mm,width=3mm]}] (x1) -- (v3);
\draw[black, arrows={->[line width=1pt,black,length=3mm,width=3mm]}] (x3) -- (y2);
\draw[black, arrows={->[line width=1pt,black,length=3mm,width=3mm]}] (y1) -- (y3);

\draw (v1) -- (v3);

\draw[anchor = north] (x1) node {{\small $x_1$}};
\draw[anchor = south] (x2) node {{\small $x_2$}};
\draw[anchor = north] (x3) node {{\small $x_3$}};

\draw[anchor = south] (y1) node {{\small $y_1$}};
\draw[anchor = north] (v3) node {{\small $w_1$}};
\draw[anchor = south] (v2) node {{\small $q$}};
\draw[anchor = north] (v1) node {{\small $s$}};

\draw[anchor = north] (y2) node {{\small $z_1$}};
\draw[anchor = north] (y3) node {{\small $r$}};
\draw (x1) [fill=red]circle (\vr);
\draw (x2) [fill=black] circle (\vr);
\draw (x3) [fill=black] circle (\vr);

\draw (v1) [fill=black]circle (\vr);
\draw (v2) [fill=red] circle (\vr);
\draw (v3) [fill=black] circle (\vr);

\draw (y2) [fill=black]circle (\vr);
\draw (y1) [fill=red] circle (\vr);
\draw (y3) [fill=white] circle (\vr);

\draw (-1,0.45) node {{\small $U_2$}};
\draw (4,0.45) node {{\small $U_1$}};
\draw (9,0.45) node {{\small $U_3$}};

\end{tikzpicture}
\end{center}
\vskip -0.5 cm
\caption{The units $U_1$, $U_2$ and $U_3$ in the proof of Claim~\ref{two}.} \label{Claim_Two}
\end{figure}
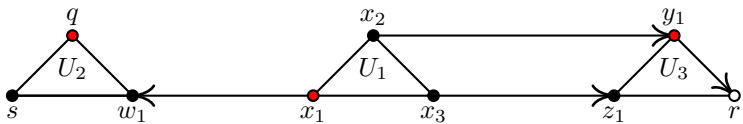

\begin{claim}
\label{two}
If $U_2 = U_4$ or $U_3 = U_4$, then $Z(G) \le \alpha(G) + 1$ and $Z(G) \le \alpha'(G)$. 
\end{claim}
\proof Renaming vertices if necessary, we may assume that $U_3 = U_4$. Thus, $\{y_1,z_1\} \subset V(U_3)$. Let $r$ be the third vertex in $U_3$. We now apply the Triangle-Rule from Claim~\ref{one} to the triangle-unit $U_2$, and greedily color a vertex, say $q \in U_2$, distinct from $w_1$, and update $S$ and $I$ by adding to these sets the vertex~$q$; that is, $S:=S\cup\{q\}$ and $I := I\cup\{q\}$. (We note that the vertices $q$ and $r$ may possibly be adjacent.) Let $s$ be the third vertex in $U_2$. This process is illustrated in Figure~\ref{Claim_Two}, where the red vertices $x_1$, $y_1$ and $q$ belong to the current independent set $I$. Recall that the vertex $x_3$ forces its $S$-uncolored neighbor $z_1$ to be colored. Thus, the vertex $y_1$ may now be played and force its uncolored neighbor $r$ in $U_3$ to be colored. By the Triangle-Rule, the vertex $q$ is played and forces one new vertex to be colored. Thus, we have ensured that each vertex in $I$ is independent and also played during the forcing process on $G$, again see Figure~\ref{Claim_Two}.

We now allow the zero forcing process to continue. At each halting point (before all of $V(G)$ becomes colored) there must be a colored vertex, say $v$, with exactly two uncolored neighbors that belong to the same triangle-unit as $v$. We now apply the Triangle-Rule which assures that our greedy coloring process will continue and eventually result in all of $V(G)$ becoming colored. Moreover, the Triangle-Rule ensures that at each halting point of the forcing process, we greedily color a vertex which is independent from any previously colored vertices. Thus, $Z(T) \le S| = |I| + 1 \le \alpha(G) + 1$. 
As in the proof of Claim~\ref{one}, every vertex in the constructed set $S$ is an $S$-forcing vertex, and the edges along which the vertices in $S$ force are independent of each other, implying that the graph $G$ contains a matching of size~$|S|$. Hence, $Z(T) \le |S| \le \alpha'(G)$.  This completes the proof of Claim~\ref{two}.~\smallqed

\begin{figure}[htb]
\begin{center}
\begin{tikzpicture}[scale=.8,style=thick,x=1cm,y=1cm, =>stealth]
\def\vr{2.5pt} 

\path (3,0) coordinate (x1);
\path (4,1) coordinate (x2);
\path (5,0) coordinate (x3);

\path (-2,0) coordinate (v1);
\path (-1,1) coordinate (v2);
\path (0,0) coordinate (v3);

\path (8,0) coordinate (y2);
\path (9,1) coordinate (y1);
\path (10,0) coordinate (y3);

\path (3,3) coordinate (w1);
\path (4,2) coordinate (w2);
\path (5,3) coordinate (w3);

\path (9.5,2) coordinate (e1);
\path (-2,2) coordinate (e2);
\path (2,4) coordinate (e3);
\path (-3,0) coordinate (e4);


\draw (x1) -- (x2);
\draw (x1) -- (x3);
\draw (x2) -- (x3);

\draw (v1) -- (v2);
\draw (v1) -- (v3);
\draw (v2) -- (v3);

\draw(y1)--(y2);
\draw(y3)--(y2);
\draw(y3)--(y1);

\draw (w1) -- (w2);
\draw (w1) -- (w3);
\draw (w2) -- (w3);

\draw[black, arrows={->[line width=1pt,black,length=3mm,width=3mm]}] (x2) -- (w2);
\draw[black, arrows={->[line width=1pt,black,length=3mm,width=3mm]}] (x1) -- (v3);
\draw[black, arrows={->[line width=1pt,black,length=3mm,width=3mm]}] (x3) -- (y2);

\draw[anchor = north] (x1) node {{\small $x_1$}};
\draw[anchor = east] (x2) node {{\small $x_2$}};
\draw[anchor = north] (x3) node {{\small $x_3$}};

\draw[anchor = west] (w2) node {{\small $y_1$}};
\draw[anchor = north] (v3) node {{\small $w_1$}};
\draw[anchor = north] (y2) node {{\small $z_1$}};

\draw (x1) [fill=red]circle (\vr);
\draw (x2) [fill=black] circle (\vr);
\draw (x3) [fill=black] circle (\vr);

\draw (v1) [fill=white]circle (\vr);
\draw (v2) [fill=white] circle (\vr);
\draw (v3) [fill=black] circle (\vr);

\draw (y2) [fill=black]circle (\vr);
\draw (y1) [fill=white] circle (\vr);
\draw (y3) [fill=white] circle (\vr);

\draw (w1) [fill=white]circle (\vr);
\draw (w2) [fill=red] circle (\vr);
\draw (w3) [fill=white] circle (\vr);

\draw (4,2.625) node {{\small $U_4$}};
\draw (-1,0.45) node {{\small $U_2$}};
\draw (4,0.45) node {{\small $U_1$}};
\draw (9,0.45) node {{\small $U_3$}};

\end{tikzpicture}
\end{center}
\vskip -0.5 cm
\caption{The distinct units $U_1$, $U_2$, $U_3$ and $U_4$.} \label{Claim_Three}
\end{figure}
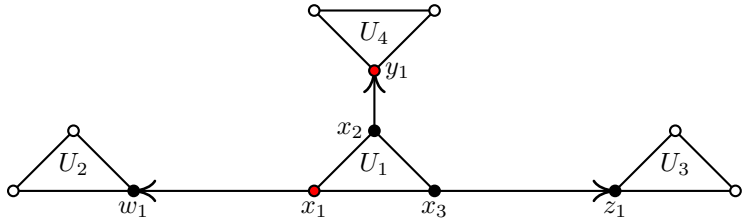

\medskip
By our earlier assumptions, every unit in $G$ is a triangle-unit. If two triangle-units in the $\Delta$-D-partition are joining by two edges, then we can choose the units $U_1$, $U_2$, and $U_3$ so that $U_2 = U_4$ or $U_3 = U_4$, and the desired result follows from Claim~\ref{two}. Hence, we may assume that every two adjacent units are joined by exactly one edge. In particular, we note that the units $U_1$, $U_2$, $U_3$ and $U_4$ are all distinct, as illustrated in Figure~\ref{Claim_Three} where the red vertices $x_1$ and $y_1$ belong to the current independent set $I$.

\vskip -0.5 cm 

We define the contraction multigraph of $G$, denoted $M_G$, to be the multigraph whose vertices correspond to the triangle-units in $G$ and where two vertices in $M_G$ are joined by the number of edges joining the corresponding triangle-units in $G$. By our earlier assumptions, adjacent triangle-units in $G$ are joined by exactly one edge. Thus, $M_G$ has no multiple edges. Further since $G$ has no diamond-unit, we note that $M_G$ has no loops. Therefore, $M_G$ is a (connected) cubic graph.

Since $M_G$ is a cubic graph, it contains at least one cycle. Let $C \colon b_0b_1 \ldots b_kb_0$ be a shortest cycle in $M_G$. For $i \in [k] \cup \{0\}$, let $T_i$ be the triangle-unit in $G$ associated with the vertex $b_i$ in $M_G$. We now choose the initially-forcing unit $U_1$ to be the triangle-unit $T_0$ in $G$ (associated with the vertex $b_0$). Renaming the units adjacent to the unit $U_1$ if necessary, we may assume that $U_2 = T_1$ and $U_4 = T_k$. Thus, $T_1, \dots, T_k$ is a sequence of distinct triangle-units where $T_i$ and $T_{i+1}$ are adjacent units for $i \in [k-1]$. Further, we note that this sequence $T_1, \dots, T_k$ of triangle-units does not contain the unit $U_1$.  

\begin{figure}[htb]
\begin{center}
\begin{tikzpicture}[scale=.8,style=thick,x=1cm,y=1cm, =>stealth]
\def\vr{2.5pt} 

\path (3,0) coordinate (x1);
\path (4,1) coordinate (x2);
\path (5,0) coordinate (x3);

\path (-2,0) coordinate (v1);
\path (-1,1) coordinate (v2);
\path (0,0) coordinate (v3);

\path (8,0) coordinate (y2);
\path (9,1) coordinate (y1);
\path (10,0) coordinate (y3);

\path (3,3) coordinate (w1);
\path (2.9,2.9) coordinate (w1p);
\path (4,2) coordinate (w2);
\path (5,3) coordinate (w3);

\path (9.5,2) coordinate (e1);
\path (-2,2) coordinate (e2);
\path (2,4) coordinate (e3);
\path (-3,0) coordinate (e4);


\draw (x1) -- (x2);
\draw (x1) -- (x3);
\draw (x2) -- (x3);

\draw (v1) -- (v2);
\draw (v1) -- (v3);
\draw (v2) -- (v3);

\draw(y1)--(y2);
\draw(y3)--(y2);
\draw(y3)--(y1);

\draw (w1) -- (w2);
\draw (w1) -- (w3);
\draw (w2) -- (w3);

\draw[black, arrows={->[line width=1pt,black,length=3mm,width=3mm]}] (x2) -- (w2);
\draw[black, arrows={->[line width=1pt,black,length=3mm,width=3mm]}] (x1) -- (v3);
\draw[black, arrows={->[line width=1pt,black,length=3mm,width=3mm]}] (x3) -- (y2);
\draw[black, arrows={->[line width=1pt,black,length=3mm,width=3mm]}] (v3) -- (v1);
\draw[black, arrows={->[line width=1pt,black,length=3mm,width=3mm]}] (w2) -- (w3);
\draw[black, arrows={->[line width=1pt,black,length=3mm,width=3mm]}] (y2) -- (y3);
\draw[black, arrows={->[line width=1pt,black,length=3mm,width=3mm]}] (y1) -- (e1);
\draw[black, arrows={->[line width=1pt,black,length=3mm,width=3mm]}] (v2) -- (e2);
\draw[black, arrows={->[line width=1pt,black,length=3mm,width=3mm]}] (e3) -- (w1);
\draw[black, arrows={->[line width=1pt,black,length=3mm,width=3mm]}] (v1) -- (e4);
\draw[anchor = north] (x1) node {{\small $x_1$}};
\draw[anchor = east] (x2) node {{\small $x_2$}};
\draw[anchor = north] (x3) node {{\small $x_3$}};

\draw[anchor = east] (w2) node {{\small $y_1 = v_k$}};
\draw[anchor = east] (w1p) node {{\small $w_k$}};
\draw[anchor = west] (w3) node {{\small $u_k$}};

\draw[anchor = north] (y2) node {{\small $z_1$}};

\draw[anchor = north] (v3) node {{\small $w_1$}};
\draw[anchor = south] (v2) node {{\small $v_1$}};
\draw[anchor = south east] (e3) node {{\small $\dots\rightarrow T_{k-1}$}};
\draw[anchor = south east] (e2) node {{\small $\dots\leftarrow T_2$}};

\draw[anchor = north] (v1) node {{\small $u_1$}}; 

\draw (x1) [fill=red]circle (\vr);
\draw (x2) [fill=black] circle (\vr);
\draw (x3) [fill=black] circle (\vr);

\draw (v1) [fill=black]circle (\vr);
\draw (v2) [fill=red] circle (\vr);
\draw (v3) [fill=black] circle (\vr);

\draw (y2) [fill=black]circle (\vr);
\draw (y1) [fill=red] circle (\vr);
\draw (y3) [fill=black] circle (\vr);

\draw (w1) [fill=black]circle (\vr);
\draw (w2) [fill=red] circle (\vr);
\draw (w3) [fill=white] circle (\vr);

\draw (-1,-0.75) node {{\small $T_1 = U_2$}}; 
\draw (4,0.45) node {{\small $U_1$}};;
\draw (9,0.45) node {{\small $U_3$}};;
\draw (4,3.55) node {{\small $T_k = U_4$}};

\end{tikzpicture}
\end{center}
\vskip -0.5 cm
\caption{A triangle-chain starting at $U_2$ and ending at $U_4$} \label{Claim_Three_One}
\end{figure}
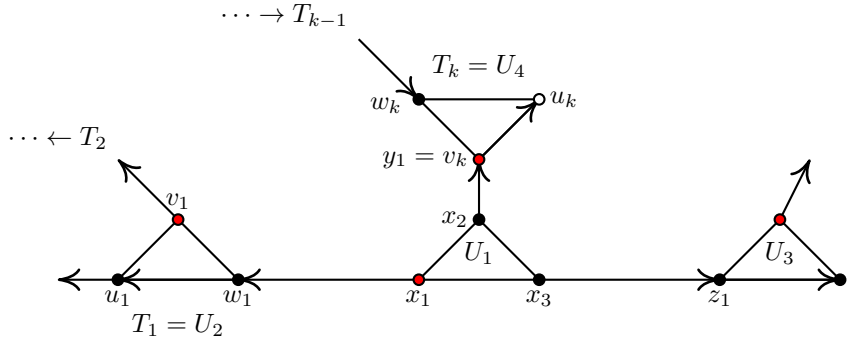

Let $V(T_i) = \{u_i,v_i,w_i\}$ for $i \in [k]$, where $w_1v_1w_2v_2 \ldots w_k$ is a $(w_1,w_k)$-path in $G$ and where $y_1 = v_k$. We now greedily color the vertex $v_i$ from the triangle-unit $T_i$ for each $i \in [k-1]$. We update $S$ and $I$ by adding to these sets the vertices $v_1,\ldots,v_{k-1}$. By the Triangle-Rule, the vertex $w_1$ is played and forces the vertex $u_1$ to be colored, and next the vertex $v_1$ is played and forces the vertex $w_2$ to be colored. Thereafter, by the Triangle-Rule, the vertex $w_2$ is played and forces the vertex $u_2$ to be colored, and next the vertex $v_2$ is played and forces the vertex $w_3$ to be colored. Continuing in this way, by the Triangle-Rule, once the vertex $w_i$ is colored, it is played and forces the vertex $u_i$ to be colored, and next the vertex $v_i$ is played and forces the vertex $w_{i+1}$ to be colored for each $i \in [k-1]$. Once the vertex $v_{k-1}$ is played and forces the vertex $w_k \in V(U_4)$ to be colored, we note that at this point of the forcing process the vertex $y_1 = v_k$ has exactly one uncolored neighbor, namely the vertex $u_k$. The vertex $y_1$ is now played and forces the vertex $u_k$ to be colored, as illustrated in  Figure~\ref{Claim_Three_One}, where the red vertices belong to the current independent set $I$.

As before, we now allow the forcing process to continue, and at each halting step (before all of $V(G)$ becomes colored) we apply the Triangle-Rule, to greedily color new vertices until all of $V(G)$ becoming colored.  Moreover, the Triangle-Rule ensures that at each halting point of the forcing process, we greedily color a vertex which is independent from any previously colored vertices. Adding the greedily colored vertices to $S$ and $I$, the resulting set $S$ becomes a zero forcing set and the resulting set $I$ is an independent set, where each vertex in $I$ is played during the forcing process. Thus, $Z(T) \le S| = |I| + 1 \le \alpha(G) + 1$.
As in the proof of Claim~\ref{one}, every vertex in the constructed set $S$ is an $S$-forcing vertex, and the edges along which the vertices in $S$ force are independent of each other, implying that the graph $G$ contains a matching of size~$|S|$. Hence, $Z(T) \le |S| \le \alpha'(G)$.  This completes the proof of Theorem~\ref{main_thm}.~\qed

\section{Closing Remarks}
\label{closing}

In this paper we have shown in Theorem~\ref{main_thm} that the zero forcing number of a connected, claw-free, cubic graph different from $K_4$ is at most its independence number plus one. However, it remains an open problem to characterize those graphs achieving equality in the upper bounds of Theorem~\ref{main_thm}. We believe that equality holds for only a finite set of connected, claw-free, cubic graphs. Indeed, we were unable to find any such graphs $G$ different form $C_3 \cp K_2$, $N_2$, and $N_3$ satisfying $Z(G) = \alpha(G) + 1$. If no such graphs exist, then this would imply that every connected, claw-free, cubic graph $G$ of order~$n \ge 14$ satisfies $Z(G) \le \alpha(G)$.

\newpage
We remark that our proof of Theorem~\ref{main_thm} shows that we can construct a zero forcing set in a connected, claw-free, cubic graph starting with three vertices from one unit and at most one vertex from every other unit. Thus, as an immediate consequence of our proof of Theorem~\ref{main_thm} we have the following upper bound on the zero forcing number of a claw-free, cubic graph of order~$n$ with $n_3$ triangle-units and $n_4$ diamond-units, noting that $n = 3n_3 + 4n_4 \ge 3(n_3 + n_4)$, and so $n_3 + n_4 \le n/3$.

\begin{cor}
\label{c:cor1} 
If $G \ne K_4$ is a connected, claw-free, cubic graph of order~$n$ with $n_3$ triangle-units and $n_4$ diamond-units, then the following holds. \\
\indent {\rm (a)} $Z(F) \le n_3 + n_4 + 2$. \\
\indent {\rm (b)} $Z(F) \le \frac{1}{3}n + 2$. 
\end{cor}

We note that if $G$ is a graph in the statement of  Corollary~\ref{c:cor1} satisfying $Z(F) = \frac{1}{3}n + 2$, then every unit in $G$ is a triangle-unit. Further, every two adjacent triangle-units in $G$ are joined by exactly one edge. However, it remains an open problem to characterize the graphs achieving equality in the upper bound of Corollary~\ref{c:cor1}(b).

\end{document}